\documentclass[10pt,reqno,a4paper]{amsart}

\usepackage{amsthm}
\usepackage{amsmath,amsfonts,amssymb}
\usepackage{graphics,color}
\usepackage{enumerate}
\usepackage{float}

\newtheoremstyle{mytheorem}{4pt}{4pt}{\itshape}{}{\scshape}{.}{.5em}{}
\theoremstyle{mytheorem}
\newtheorem{theorem}{Theorem}

\newtheorem{corollary}{Corollary}
\newtheorem{definition}{Definition}

\newtheorem{proposition}{Proposition}

\newlength{\barwidth}
\newlength{\intwidth}
\newlength{\mylen}

\newcommand*{\eps}{\varepsilon}
\newcommand*{\R}{\ensuremath{\mathbb{R}}}

\renewcommand*{\div}{\ensuremath{\textrm{div }}}

\renewcommand*{\P}{\ensuremath{\mathcal{P}}}
\renewcommand*{\S}{\ensuremath{\mathcal{S}}}
\newcommand*{\supp}{\ensuremath{\mathrm{supp\,}}}

\newcommand*{\e}{\ensuremath{\varepsilon}}

\begin{document}

\title{Weak-strong uniqueness for measure-valued solutions}

\author{Yann Brenier}
\address{CNRS\\ 
FR 2800 Wolfgang D\"oblin\\
Universit\'e de Nice\\
France}
\email{brenier@math.unice.fr}

\author{Camillo De Lellis}
\address{Institut f\"ur Mathematik\\
         Universit\"at Z\"urich\\
         CH-8057 Z\"urich\\
         Switzerland}
\email{camillo.delellis@math.unizh.ch}

\author{L\'aszl\'o Sz\'ekelyhidi Jr.}
\address{Hausdorff Center for Mathematics\\
         Universit\"at Bonn\\
         D-53115 Bonn\\
         Germany}
\email{laszlo.szekelyhidi@hcm.uni-bonn.de}


\maketitle

\begin{abstract}
We prove the weak-strong uniqueness for measure-valued solutions of the incompressible Euler equations. These were introduced by R.DiPerna and A.Majda in their landmark paper \cite{DiPernaMajda}, where in particular global existence to any $L^2$ initial data was proven. Whether measure-valued solutions agree with classical solutions if the latter exist has apparently remained open. 

We also show that DiPerna's measure-valued solutions to systems of conservation laws 
have the weak-strong uniqueness property.
\end{abstract}

\maketitle

\section{Introduction}
In \cite{DiPerna:mvsolutions} DiPerna introduced the notion of measure-valued solutions to 
conservation laws, following the pioneering work of L. Tartar on compensated compactness and Young measures. DiPerna 
worked in the context of $L^{\infty}$ solutions and thus probabilities in state space which 
are compactly supported. 
While this is sufficient in one space dimension, in general one only has a uniform energy bound, usually $L^2$, to 
work with. This is the case in particular for the incompressible Euler equations. In \cite{DiPernaMajda} 
DiPerna and Majda extended the notion of measure-valued solutions to this unbounded case. In \cite{LionsBook}
Lions remarked that for any reasonable notion of generalized solution one should require
a weak-strong uniqueness property: any time that the Cauchy problem has a ``classical'' solution, the
generalized ones should coincide with it. Lions observed that such a result is not known for 
the DiPerna-Majda's solutions and he introduced his ``dissipative solutions'', for which he could prove existence and 
weak-strong uniqueness. The remark that the weak-strong uniqueness does not seem to hold in
the DiPerna-Majda's framework
has been taken up by several other authors in the literature (see for instance 
\cite{Necas}).

Since the pioneering work of Scheffer \cite{Scheffer}, 
it is well--known that not even distributional solutions
to the Euler equations satisfy Lions' weak-strong uniqueness requirement (see also \cite{Shnirelman}, 
\cite{DeSz1} and \cite{DeSz2}). It is therefore necessary to introduce some form of energy conservation
in order to hope for this property. We show in this paper that this can be done successfully even along the 
ideas of DiPerna and Majda. Namely, it is possible to introduce a notion
of ``admissible measure-valued solution'' for which existence and weak-strong uniqueness holds.
In fact our argument shows that the barycenters of such solutions (see below for the relevant definitions)
are dissipative solutions in the sense of Lions (note, however, that the ultimate conclusion of the proof 
is that the entire measure-valued solution, and not only its barycenter, coincides with the classical one).
An interesting corollary of this analysis is that, whenever the Cauchy problem for the Euler equations
has a solution with a certain minimum regularity (slightly weaker than Lipschitz), any
sequence of Leray solutions to the vanishing viscosity approximation must converge to it. 
Known results in the literature about the convergence of solutions of Navier-Stokes to Euler (see for instance \cite{Constantin,Masmoudi}) assume more regularity. 

This paper has been inspired by the works of Brenier and Grenier \cite{Brenier,BrenierGrenier}. The main idea
of the arguments is taken from these papers and
it is a modification of a classical energy method which works for a variety
of systems of evolutionary partial differential
equations in conservation form. Our contribution is essentially 
of technical nature, clarifying the correct functional-analytic
framework to make this idea work: note, indeed, that, besides the introduction of a suitable 
energy inequality, our definition of measure-valued solutions has some other substantial differences
from the one of DiPerna and Majda. We conclude the note by showing that the same remark
can be easily extended to hyperbolic systems of conservation laws which have a strictly convex entropy. 
Namely, the well--known works of Dafermos and DiPerna (see for instance \cite{DafermosBook} Theorem 5.3.1)
can be extended to DiPerna's measure valued solutions, once we assume a suitable entropy condition. 
The proof of this statement is contained in Section \ref{s:CL}.
The other two Sections \ref{s:GYM} and \ref{s:Euler} discuss, respectively,
generalized Young measures and the results mentioned above for incompressible Euler.

\section{Generalized Young measures}\label{s:GYM}

Let $\Omega\subset\R^m$ be an open set and consider a bounded sequence $\{u_j\}\subset L^p(\Omega;\R^n)$. 
DiPerna and Majda defined generalized Young measures in order to describe weak limits of the form
$$
\lim_{j\to\infty}\int_{\Omega}\phi(y)g\bigl(u_j(y)\bigr)\,dy
$$
with $\phi\in C_c (\Omega)$ and the test function $g$ is of the form
\begin{equation}\label{e:bc}
g(\xi)=\tilde g(\xi)(1+|\xi|^p)\quad\textrm{ for some }\tilde g\in BC(\R^n).
\end{equation}
Here $BC(\R^n)$ denotes the set of bounded continuous functions on $\R^n$, and so \eqref{e:bc} defines the largest class of test functions for which one expects to be able to represent the weak limit of $g(u_j)$. 
Since $BC(\R^n)$ is isometrically isomorphic to $C(\beta\R^n)$, where $\beta\R^n$ is the Stone-\v Cech compactification of $\R^n$, the most general way to represent the weak limits is using a measure $\hat\nu$ in the space  
$$
\mathcal{M}(\Omega\times \beta\R^n)=C_0(\Omega\times\beta\R^n)^*.
$$
In other words, there exists a subsequence $u_{j_k}$ such that
$$
\lim_{k\to\infty}\int_{\Omega}\phi(y)\tilde{g}\bigl(u_{j_k}(y)\bigr)(1+|u_j(y)|^p)dy\,=\,\langle \hat\nu,\phi\tilde{g}\rangle 
$$
for all $\phi\in C_0(\Omega)$ and $\tilde g\in BC(\R^n)$. The measure $\hat\nu$ is called the \emph{generalized Young measure}, c.f. \cite[Corollary 4.1]{DiPernaMajda}.

Moreover, DiPerna and Majda proved in \cite[Theorem 4.3]{DiPernaMajda} that for a certain subclass of test functions the measure $\hat\nu$ admits 
a disintegration into a family of probability measures. More precisely, 
let $\mathcal{F}\subset BC(\R^n)$ be a separable completely regular subalgebra, and 
let $\sigma\in\mathcal{M}(\Omega)$ the projection onto $\Omega$ of $\hat\nu$, i.e.
$$
\sigma(E)=\hat\nu(E\times\beta\R^n)\textrm{ for }E\subset\Omega.
$$
There exists a $\sigma$-measurable map 
$$
\Omega\to{\rm Prob}(\beta_{\mathcal{F}}\R^n):\,y\mapsto\hat\nu_y
$$
such that for every $\tilde{g}\in\mathcal{F}$ and $\phi\in C_0(\Omega)$
$$
\langle \hat\nu,\phi\tilde g\rangle = \int_{\Omega}\phi\int_{\beta_{\mathcal{F}}\R^n}\tilde g\,d\hat\nu_y\,d\sigma.
$$
A particularly useful class of test functions is 
$$
\mathcal{F}=\{\tilde g\in BC(\R^n):\,\tilde g^{\infty}(\xi):=\lim_{s\to\infty}\tilde g(s\xi)\textrm{ exists and is continuous on }S^{n-1}\}
$$
In this case $\beta_\mathcal{F}\R^n$ can be identified with the closed unit ball $\overline{B^n}$. Observe that in this case
$g^{\infty}(\xi)=\tilde g^{\infty}(\xi)(1+|\xi|^p)$ coincides with the $L^p$-recession function
$$
g^{\infty}(\theta)=\lim_{s\to\infty}\frac{g(s\theta)}{s^p}\,\textrm{ for all }\theta\in S^{n-1}.
$$
A further step in the analysis of such measures was taken by Alibert and Bouttich\'e in \cite{AlibertBouchitte}. They obtain a decomposition of $\nu$ into an triple 
$$
(\nu_y,\nu_y^\infty,\lambda)\in {\rm Prob}(\R^n)\times {\rm Prob}(S^{n-1})\times \mathcal{M}^+(\Omega),
$$
such that 
$$
\int_{\Omega}\phi\int_{\beta_{\mathcal{F}}\R^n}\tilde g\,d\hat\nu_y\,d\sigma = \int_{\Omega}\phi\int_{\R^n}g\,d\nu_y\,dy + \int_{\Omega}\phi\int_{S^{n-1}}g^{\infty}\,d\nu_y^\infty\,d\lambda\quad\textrm{ for all }\tilde g\in\mathcal{F}.
$$
This is obtained by using the observation (by testing with $g(\xi)\equiv 1$) that 
$$
\textrm{ for }\sigma_s-{\rm a.e. }\,y:\quad\hat\nu_y\textrm{ is concentrated on }\beta_\mathcal{F}\R^n\setminus\R^n,
$$
where $\sigma_s$ is the singular part of $\sigma$ with respect to Lebesgue measure. After appropriate normalizations one is lead to a representation of the above form. 

At this point we introduce the following notation, adapted from \cite{AlibertBouchitte}: given a Radon measure $\lambda$ and a topological space $X$, we denote by $\P(\lambda;X)$ the set of parametrized families of probability measures $(\nu_y)$ on $X$ which depend $\lambda$-measurably on the parameter $y$. In the particular case when $\lambda$ is Lebesgue measure on $\Omega\subset\R^m$, we write $\P(\Omega;X)$.

In summary, one has the following result:
\begin{theorem}[DiPerna-Majda, Alibert-Bouchitt\'e]
Let $\{u_k\}$ be a bounded sequence in $L^p(\Omega;\R^n)$. There exists a subsequence $\{u_{k_j}\}$, a 
nonnegative Radon measure $\lambda$ and parametrized families of probability measures $\nu\in\P(\Omega;\R^n)$,
$\nu^{\infty}\in\P(\lambda;\S^{n-1})$ such that:
\begin{equation}\label{e:generation}
g(u_{k_j})\,\overset{*}\rightharpoonup\,\langle\nu,g\rangle +\langle\nu^{\infty},g^{\infty}\rangle \lambda
\end{equation}
in the sense of measures, for every $g\in\mathcal{F}$.
\end{theorem}

\section{Admissible measure-valued solutions of Euler}\label{s:Euler}

Let $v_0\in L^2(\R^n)$ with $\div v_0=0$. 
Following \cite{DiPernaMajda}, we consider a sequence of Leray solutions $v_{\e}\in L^{\infty}(\R_+;L^2(\R^n))$ with vanishing viscosity $\e\to 0$. Using the uniform energy bound
$$
\int_{\R^n}|v_\e(x,t)|^2\,dx\leq \int_{\R^n}|v_0(x)|^2\,dx
$$ 
it is easy to see that for any bounded $\Omega\subset\R_+\times\R^n$ 
a suitable subsequence generates a measure-valued solution. Then, by considering a standard diagonal argument we can extend this to all of $\R_+\times\R^n$. Using the representation above for the generalized Young measure $\hat\nu$, we thus obtain  a triple $(\nu,\nu^{\infty},\lambda)$ with
$\lambda\in\mathcal{M}^+(\R_+\times\R^n)$ and 
$$
\nu\in\P(\R_+\times\R^n;\R^n),\,\nu^{\infty}\in\P(\lambda;\S^{n-1}),
$$
such that the equations 
\begin{eqnarray}
\partial_t\langle \nu,\xi\rangle + \div\Bigl(\langle\nu,\xi\otimes\xi\rangle+\langle\nu^{\infty},\theta\otimes\theta\rangle\,\lambda\Bigr)+\nabla p&=&0,\label{e:euler1}\\
\div\langle\nu,\xi\rangle&=&0\label{e:euler2}
\end{eqnarray}
hold in the sense of distributions. Here the bracket $\langle\cdot,\cdot\rangle$ denotes the appropriate integrals, so that
$$
\langle \nu,\xi\otimes\xi\rangle = \int_{\R^n}(\xi\otimes\xi)\,\nu_{x,t}(d\xi)\,,\quad
\langle \nu^{\infty},\theta\otimes\theta\rangle = \int_{S^{n-1}}(\theta\otimes\theta)\,\nu_{x,t}^{\infty}(d\theta),
$$
and in particular
$$
\bar{\nu}(x,t) :=\langle\nu_{x,t},\xi\rangle 
$$
stands for the barycenter of the probability measure $\nu_{x,t}$.

Now, testing \eqref{e:generation} with $g(\xi)=|\xi|^2$ (and hence $g^{\infty}(\theta)\equiv 1$) and using the energy bound for
the Leray solutions $v_\varepsilon$, we obtain
\begin{equation}\label{e:energy}
\begin{split}
\int_{\R_+}\int_{\R^n}\varphi(x)\chi(t)\langle\nu_{x,t},|\xi|^2\rangle \,dxdt
+\int_{\R_+}&\int_{\R^n}\varphi(x)\chi(t)\, d\lambda\leq\\
& \leq\|\varphi\|_{\infty}\|\chi\|_{1}\int_{\R^n}|v_0(x)|^2\,dx
\end{split}
\end{equation}
for all $\varphi\in C_c(\R^n)$ and $\chi\in C_c(\R_+)$.
Jensen's inequality and the first term implies that $\bar{\nu}\in L^{\infty}_tL^2_x$, whereas the second term and a standard slicing argument implies that $\lambda$ admits the representation
$$
\lambda=\lambda_t(dx)\otimes dt,
$$
where $t\mapsto \lambda_t$ is a measurable $\mathcal{M}_+(\R^n)$-valued function. Thus, we may define the energy of the generalized Young measure as $E\in L^{\infty}(\R_+)$ by 
$$
E(t)=\frac{1}{2}\int_{\R^n}\langle\nu_{x,t},|\xi|^2\rangle dx+\frac{1}{2}\lambda_t(\R^n)
$$
and obviously from \eqref{e:energy} we conclude
\begin{equation}\label{e:energy2}
E (t) \;\leq\; \frac{1}{2} \int_{\R^n} |v_0 (x)|^2\,dx
\qquad \mbox{for a.e. $t$.}
\end{equation} 
Moreover, 
from \eqref{e:euler1} we deduce that $\bar{\nu}$ can be redefined on a set of times of measure zero so that for any $\varphi\in L^2(\R^n)$ the function
$$
t\mapsto \int_{\R^n}\varphi(x)\bar{\nu}(x,t)\,dx
$$
is continuous. Hence we may assume that $\bar{\nu}\in C([0,\infty[;L^2_w(\R^n))$ and in particular $\bar{\nu}(\cdot,t)\rightharpoonup v_0(\cdot)$ in $L^2$ as $t\to 0$.
We can combine this information with \eqref{e:euler1} in the form
\begin{equation}\label{e:euler3}
\begin{split}
\int_0^\infty\int_{\R^n}\partial_t\phi\cdot\bar\nu+ \nabla\phi:\langle\nu,\xi\otimes\xi\rangle\,dxdt
+&\int_0^\infty\int_{\R^n}\nabla\phi:\langle\nu^{\infty},\theta\otimes\theta\rangle\,\lambda(dx,dt)
\\
&\quad=-\int_{\R^n}\phi(x,0)v_0(x)\,dx
\end{split}
\end{equation}
for all $\phi\in C^{\infty}_c([0,\infty[\times\R^n;\R^n)$ with 
$\div\phi=0$ (we use here the common
notation $A:B = \sum_{ij} A_{ij} B_{ij}$).

Motivated by the above, in analogy with DiPerna \cite[Section 4b)]{DiPerna:mvsolutions} we make the following definition:
\begin{definition}
A triple $(\nu,\nu^\infty,\lambda)$
is an admissible measure-valued solution of the Euler equations with initial data $v_0$ provided \eqref{e:euler2},\eqref{e:energy} and \eqref{e:euler3} hold.
\end{definition}
In the above we have shown, in particular, the following

\begin{proposition}\label{p:existence}
For any initial data $v_0\in L^2(\R^n)$,
any sequence of Leray's solutions to Navier-Stokes with vanishing viscosity has
a subsequence converging to an admissible measure-valued solution. There exists, therefore,
at least one such solution. 
\end{proposition}

\subsection{Weak-strong uniqueness}

Let $v_0\in L^2(\R^n)$ with $\div v_0=0$, and consider the initial value problem for the incompressible Euler equations. We show here the following theorem.

\begin{theorem}\label{t:main}
Assume that $v\in C([0,T];L^2(\R^n))$ is a solution with 
\begin{equation}\label{e:symmetric}
\int_0^T\|\nabla v+\nabla v^T\|_{\infty}\,dt<\infty
\end{equation}
and let $(\nu,\nu^\infty,\lambda)$ be any admissible measure-valued solution. 
Then $\lambda=0$ and $\nu_{x,t}=\delta_{v(x,t)}$ for a.e. $(x,t)$.
\end{theorem}

Indeed, the proof below yields easily the following proposition.

\begin{proposition}\label{p:dissipative}
Let $(\nu,\nu^\infty,\lambda)$ be an admissible measure-valued solution. Then
$v (t,x):= \langle \nu_{t,x}, \xi\rangle$ is a dissipative solution in the sense of
Lions. 
\end{proposition}

Finally, we observe that Proposition \ref{p:existence} and Theorem \ref{t:main}
has the following interesting corollary.

\begin{corollary}\label{c:convergence}
Assume that, for some divergence-free $v_0\in L^2$ 
there is a solution $v\in C([0,T];L^2(\R^n))$ of Euler such that
\eqref{e:symmetric} holds.
Then, any sequence of Leray's solutions to the corresponding vanishing 
viscosity approximation converge to $v$ in $L^2((0,T)\times \R^n)$.
\end{corollary}

\begin{proof}[Proof of Theorem \ref{t:main}]
Let 
$$
F(t)=\frac{1}{2}\int_{\R^n}\langle \nu_{x,t},|\xi-v(x,t)|^2\rangle\,dx+\frac{1}{2}\lambda_t(\R^n).
$$
Furthermore, for $\varphi\in C_c^\infty(\R^n)$ with $\div\varphi=0$ define
$$
F^{\varphi}(t)=\frac{1}{2}\int_{\R^n}\varphi(x)\langle \nu_{x,t},|\xi-v(x,t)|^2\rangle\,dx+\frac{1}{2}\int_{\R^n}\varphi(x)\,\lambda_t(dx).
$$
Observe that $F^{\varphi}\in L^{\infty}(0,T)$ by \eqref{e:energy}. Let $\chi\in C^{\infty}_c (0,T)$ and 
consider
\begin{equation}\label{e:F}
\begin{split}
\int_0^T\chi'(t)\,F^{\varphi}(t)\,dt=&\frac{1}{2}\iint\chi'\varphi\langle \nu,|\xi|^2\rangle\,dxdt+\frac{1}{2}\iint\chi'\varphi\,\lambda_t(dx)dt\\
&+\frac{1}{2}\iint\chi'\varphi|v|^2\,dxdt-\iint\chi'\varphi\,\bar\nu\cdot v\,dxdt.
\end{split}
\end{equation}
Using \eqref{e:euler3} the final term above can be written as
\begin{equation}
\begin{split}
-\iint\chi'\varphi\,\bar\nu\cdot v\,dxdt&=\iint -\partial_t(\chi\varphi v)\cdot\bar\nu - \chi\varphi\bar{\nu}\cdot\bigl(\div(v\otimes v)+\nabla p\bigr)\,dxdt\\
=&\iint \chi\nabla(\varphi v):\langle\nu,\xi\otimes\xi\rangle-\chi\varphi\bar\nu\cdot\div(v\otimes v)\,dxdt\\
+&\iint \chi\nabla(\varphi v):\langle\nu^{\infty},\theta\otimes\theta\rangle\,\lambda_t(dx)dt\\-&\iint \chi\varphi\bar\nu\cdot\nabla p\,dxdt.
\end{split}
\end{equation}
Next, we use the identities
\begin{eqnarray*}
\nabla(\varphi v):(\bar{\nu}\otimes v)&=&\varphi\bar{\nu}\cdot\div(v\otimes v)-(\nabla\varphi\cdot v)(v\cdot \bar{\nu}),\\
\nabla(\varphi v):(v\otimes\bar{\nu})&=&\nabla\varphi\cdot\bar{\nu}\frac{|v|^2}{2}+\nabla(\varphi\frac{|v|^2}{2})\cdot\bar{\nu},\\
\nabla(\varphi v):(v\otimes v)&=&\nabla\varphi\cdot v\frac{|v|^2}{2}+\nabla(\varphi\frac{|v|^2}{2})\cdot v,
\end{eqnarray*}
together with $\div v=0$, $\div\bar{\nu}=0$ to rearrange further as
\begin{equation}
\begin{split}
-\iint\chi'\varphi\,\bar\nu\cdot v&\,dxdt=\iint \chi\nabla(\varphi v):\langle\nu,(\xi-v)\otimes(\xi-v)\rangle\,dxdt\\
+&\iint \chi\nabla(\varphi v):\langle\nu^{\infty},\theta\otimes\theta\rangle\,\lambda_t(dx)dt\\
+&\iint \chi\nabla\varphi\cdot\left((\bar{\nu}-v)\frac{|v|^2}{2}+\bar{\nu}p\right)+\chi(\nabla\varphi\cdot v)(v\cdot\bar{\nu})\,dxdt.
\end{split}
\end{equation}
Observe that $\|v\|_{L^{\infty}(\R^n)}$ can be bounded in terms of $\|\nabla v+\nabla v^T\|_{L^{\infty}(\R^n)}+\|v\|_{L^2(\R^n)}$. Indeed, for any ball $B_1(x_0)$
Korn's inequality implies a bound on $\|\nabla v\|_{L^p(B_1(x_0))}$, and from here the bound on $\|v\|_{L^{\infty}(B_1(x_0))}$ follows from the Sobolev embedding and the uniform $L^2$ bound. 
In turn, from the uniform $L^{\infty}$ and $L^2$ bounds on $v$ follows that $v\in L^4(\R^n)$ and $p\in L^2(\R^n)$.    
Next, take a sequence $\{\varphi_k\}$ such that
$0\leq\varphi_k(x)\leq 1$, $\varphi_k\equiv 1$ on $B_k (0)$
and $\|\nabla\varphi_k\|_{C^0}$ is uniformly bounded. 
Using dominated convergence and the bounds obtained above we see that under the assumption \eqref{e:symmetric}
$$
\iint \chi\nabla\varphi_k\cdot\left((\bar{\nu}-v)\frac{|v|^2}{2}+\bar{\nu}p\right)+\chi(\nabla\varphi_k\cdot v)(v\cdot\bar{\nu})\,dxdt\to 0,
$$
and therefore
\begin{equation}
\begin{split}
-\iint\chi'\varphi_k\,\bar\nu\cdot v\,dxdt \overset{k\to\infty}{\longrightarrow} \iint& \chi\nabla v:\langle\nu,(\xi-v)\otimes(\xi-v)\rangle\,dxdt\\
+&\iint \chi\nabla v:\langle\nu^{\infty},\theta\otimes\theta\rangle\,\lambda_t(dx)dt
\end{split}
\end{equation}
Passing to the limit also in \eqref{e:F} and symmetrizing the $\nabla u$ terms we obtain
\begin{equation}
\begin{split}
\int_0^T\chi'(t)F(t)\,dt=&\int_0^T\chi'(t)E(t)\,dt+\frac{1}{2}\int_0^T\chi'\,\int_{\R^n}|v|^2dx\,dt\\
+& \frac{1}{2}\iint \chi(\nabla v+\nabla v^T):\langle\nu,(\xi-v)\otimes(\xi-v)\rangle\,dxdt\\
+&\frac{1}{2}\iint \chi(\nabla v+\nabla v^T):\langle\nu^{\infty},\theta\otimes\theta\rangle\,\lambda_t(dx)dt
\end{split}
\end{equation}
Since $\nabla v \in L^1 ([0,T], L^q
(B))$ for every $q<\infty$ and
$p\in L^2$, it is easy to see that $\partial_t |v|^2 + {\rm div}
[(|v|^2 +2p) v] = 0$. On the other hand, integrating this
identity in space, the bounds above imply that $\int |v|^2 (x,t) dx$
is constant. Hence we deduce
$$
-\int_0^T\chi'(t)F(t)\,dt\leq -\int_0^T\chi'(t)E(t)\,dt+C\int_0^T\chi(t)\|\nabla v(t)+\nabla v(t)^T\|_{\infty}F(t)\,dt.
$$
Therefore, for almost every $s,t\in(0,T)$ 
\begin{equation}\label{e:st}
F(t)-F(s)\leq E(t)-E(s)+C\int_s^t\|\nabla v(\tau)+\nabla v(\tau)^T\|_{\infty}F(\tau)\,d\tau.
\end{equation}
Finally, observe that 
$$
F(s)=E(s)-\int_{\R^n}\bar{\nu}\cdot v\,dx+\frac{1}{2}\int_{\R^n}|v|^2dx,
$$
so that \eqref{e:st} becomes
\begin{equation*}
\begin{split}
F(t)\leq & E(t)-\int_{\R^n}\bar{\nu}(x,s)\cdot v(x,s)\,dx+
\frac{1}{2}\int_{\R^n}|v(x,s)|^2\,dx\\
&+C\int_s^t\|\nabla v(\tau)+\nabla v(\tau)^T
\|_{\infty}F(\tau)\,d\tau\, .
\end{split}
\end{equation*}
Now passing to the limit $s\to 0$ (justified since $\bar{\nu}\in CL^2_w$) 
$$
F(t)\leq E(t)-\frac{1}{2}\int_{\R^n}|v_0(x)|^2\,dx +C\int_s^t\|\nabla v(\tau)+\nabla v(\tau)^T\|_{\infty}F(\tau)\,d\tau,
$$
from which (recalling \eqref{e:energy2})
$$
F(t)\leq C\int_0^t\|\nabla v(\tau)+\nabla v(\tau)^T\|_{\infty}F(\tau)\,d\tau
$$
follows by the admissibility assumption. Finally, this last inequality implies that $F(t)=0$ for a.e. $t$, as required.
\end{proof}

\section{Hyperbolic systems of conservation laws}\label{s:CL}

In this section we consider hyperbolic systems of conservation laws
\begin{equation}\label{e:CL}
\partial_t U + {\rm div}_x F(U) =0
\end{equation}
where $U: \Omega\subset\R\times \R^n \to \R^k$ is the unknown vector
function and $F: \R^k\to \R^{n\times k}$ a $C^2$ map. Equation \eqref{e:CL}
reads therefore
$$
\partial_t u^i + \partial_{x_j} (F^{ij} (u)) \;=\; 0\, ,
$$
which for differentiable solutions becomes 
$\partial_t u^j + \partial_l F^{ij} (u) \partial_{x_j} u^l = 0$ (in these
last identities and in what follows we use Einstein's summation
convention on repeated indices).

We assume that \eqref{e:CL} has a strictly convex entropy, i.e. that there
is a $C^2$ map $(\eta,q):\R^k\to \R\times \R^n$ 
such that $D^2 \eta\geq c_0 {\rm Id}>0$ and 
\begin{equation}\label{e:def_entropy}
\partial_i \eta \partial_l F^{ij} = \partial_l q^j\, .
\end{equation}
Thus, any Lipschitz solution of \eqref{e:CL} satisfies the identity
$\partial_t (\eta (u)) + {\rm div}_x (q(u))=0$.

\begin{definition} A bounded admissible measure-valued solution $\nu$ of \eqref{e:CL} with initial
data $U_0\in L^\infty$ is a parametrized family of propability measures
$\nu \in \P ([0,T]\times \R^n; \R^k)$ such that
\begin{itemize}
\item $t \mapsto \langle \nu_{t, \cdot}, \xi\rangle$ is a weakly$^*$ continuous map,
taking values in $L^\infty (\R^n)$;
\item the identity
\begin{equation}\label{e:CLdistrib}
\left\{\begin{array}{l}
\partial_t \langle \nu, \xi\rangle + {\rm div}_x \langle
\nu, F (\xi)\rangle \;=\; 0\\ \\
\langle \nu_{0,x}, \xi\rangle = U_0 (x)
\end{array}\right.
\end{equation}
holds in the sense of distributions;
\item the inequality
\begin{equation}\label{e:Ent}
\left\{\begin{array}{l}
\partial_t \langle \nu, \eta (\xi)\rangle + {\rm div}_x
\langle \nu, q(\xi)\rangle \;\leq\; 0\\ \\
\langle \nu_{0,x}, \eta (\xi)\rangle \;=\; \eta (U_0 (x))
\end{array}\right.
\end{equation}
holds in the sense of distributions.
\end{itemize}
\end{definition}

\begin{theorem}\label{t:CL}
Assume $U:[0,T]\times \R^n \to \R^k$ is a bounded Lipschitz solution 
of \eqref{e:CL} and $\nu$ a bounded admissible measure valued solution of
$\eqref{e:CL}$ with initial data $U_0= U(0,\cdot)$. Then $\nu_{t,x}
=\delta_{U(t,x)}$ for a.e. $(t,x)\in [0,T]\times \R^n$.
\end{theorem}

The proof follows essentially the computations of pages 98-100 in
\cite{DafermosBook}.

\begin{proof}
We start by defining the following functions of $t$ and $x$:
\begin{eqnarray}
h &:=& \langle \nu, \eta (\xi)\rangle - \eta (U)
- D\eta (U)\cdot \bigl[\langle \nu, \xi\rangle - U\bigr]\\
Y^\alpha &:=& \langle \nu, q^\alpha (\xi)\rangle -
q^\alpha (U) - \partial_l \eta (U) \big[\langle \nu, F^{l\alpha} (\xi)\rangle -
F^{l\alpha} (U)\big]\label{e:Y}\\
Z^\alpha_\beta &:=& \partial_{\beta j} \eta (U)
\big[\langle \nu, F^{j\alpha} (\xi)\rangle - F^{j\alpha} (U) - 
\partial_l F^{j\alpha} (U) \bigl(\langle \nu, \xi^l\rangle - U^l\bigr) \big]
\end{eqnarray} 
Recall that $\supp (\nu_{t,x})$ and $|U(t,x)|$ are both uniformly
bounded and that $\eta$, $q$ and $F$ are $C^2$ functions. So, there exists a 
constant $C$ such that the following identities hold for every $\xi\in \supp (\nu_{t,x})$:
\begin{eqnarray*}
|q^\alpha (U(t,x)) - q^\alpha (\xi) - \partial_i q^\alpha (U(t,x))
(U^i (t,x)-\xi)|&\leq& C |U(t,x)-\xi|^2\label{e:C2_1}\\
\left|F^{j\alpha} (U(t,x)) - F^{j\alpha} (\xi) 
- \partial_i F^{j\alpha} (U(t,x)) (U^i (t,x) - \xi^i)\right| &\leq&
C |U(t,x)-\xi|^2\label{e:C2_2}
\end{eqnarray*}
(we underline that $C$ is a constant independent of $t$,$x$ and $\xi$).
 
Plugging these last identities into \eqref{e:Y} and recalling
\eqref{e:def_entropy}, we conclude 
\begin{equation}\label{e:one_side}
|Y(t,x)|\;\leq\; C \int |U(t,x)-\xi|^2 d\nu_{t,x} (\xi)\, .
\end{equation}
On the other hand, using that $D^2 \eta\geq c_0 Id$, we easily infer
\begin{equation}\label{e:other_side}
|h (t,x)|\;\geq\; \frac{c_0}{2} \int |U(t,x)-\xi|^2 d\nu_{t,x}\, 
\end{equation}
and hence that
\begin{equation}\label{e:dafermos_5.3.5}
|Y (t,x)|\;\leq\; C_0 |h(t,x)|\, .
\end{equation}
A similar computation yields 
\begin{equation}\label{e:dafermos_5.3.5_bis}
|Z (t,x)|\;\leq\; C_1 |h(t,x)|\, .
\end{equation}
Next recall that 
\begin{equation}\label{e:classical_equality}
\partial_t (\eta (U))
+ {\rm div}_x (q(U))\;=\; 0
\end{equation} 
(because $U$ is Lipschitz). 
Fix a test function $\psi\in C^\infty_c (\R^n\times ]-T, T[)$.
Combining \eqref{e:classical_equality} with \eqref{e:Ent}, 
we conclude
\begin{eqnarray}
\int_0^T\int \left[\partial_t \psi\, h + \partial_{x_\alpha} \psi\, Y^\alpha\right]\,
&\geq& - \int_0^T \int \Big[\partial_t \psi\, \partial_i \eta (U)
\bigl(U^i - \langle \nu, \xi^i\rangle\bigr)\nonumber\\ 
&& \qquad+ \partial_{x_\alpha} \psi\, \partial_i \eta (U) \bigl(F^{i\alpha} (U) -
\langle \nu, F^{i\alpha} (\xi)\rangle \Big]\, 
\label{e:dafermos_5.3.6}
\end{eqnarray}
(no boundary term appears because
the initial condition is the same for both $\langle \nu, \eta (\xi)\rangle$
and $\eta (U)$).

In fact, by an easy approximation argument, \eqref{e:dafermos_5.3.6} holds
for any test function which is just Lipschitz continuous.
Similarly, we can use the test function $\Phi:= \psi\, D\eta (U)$ (which is
Lipschitz and compactly supported) on the identity \eqref{e:CLdistrib}
to get
\begin{equation}\label{e:dafermos_5.3.7}
\int \int \Big[ \partial_t (\psi\, \partial_i \eta (U)) \bigl( U^i -
\langle \nu, \xi^i\rangle\bigr) +
\partial_{x_\alpha} (\psi\, \partial_i \eta (U)) \bigl(F^{i\alpha} (U)
- \langle \nu, F^{i\alpha} (\xi)\rangle\Big] \;=\; 0\, .
\end{equation}
Since $U$ is Lipschitz, we can use the chain rule and \eqref{e:CL} to compute
\begin{equation}\label{e:dafermos_5.3.9}
\partial_t (\partial_i \eta (U)) \bigl( U^i -
\langle \nu, \xi^i\rangle\bigr) +
\partial_{x_\alpha} (\partial_i \eta (U)) \bigl(F^{i\alpha} (U)
- \langle \nu, F^{i\alpha} (\xi)\rangle\Big]\;=\;
\partial_{x_\alpha} U^i Z^\alpha_i
\end{equation}
Combining \eqref{e:dafermos_5.3.6}, \eqref{e:dafermos_5.3.7}
and \eqref{e:dafermos_5.3.9} we infer
\begin{equation}\label{e:dafermos_5.3.10}
\int\int \bigl[ \partial_t \psi\, h + \partial_{x_\alpha}\psi\, Y^\alpha\bigr]
\;\geq\; \int\int \psi\, \partial_{x_\alpha} U^i Z^\alpha_i\, .
\end{equation}
Next, fix any point $\tau<T$, any radius $R>0$ and $\eps\in ]0, T-\tau[$. Consider 
the test function $\psi (t,x)=\omega (t) \chi (t,x)$ where 
$$
\omega (t) \;:=\;
\left\{\begin{array}{ll}
1 & \mbox{for $0\leq t<\tau-\eps$}\\
1- \eps^{-1} (t-\tau+\eps) & \mbox{for $\tau-\eps\leq t \leq \tau$}\\
0 & \mbox{for $t\geq \tau$}.
\end{array}\right.
$$
$$
\chi (x,t) \;:=\; 
\left\{\begin{array}{ll}
1 & \mbox{if $|x|\leq R+C_0 (\tau-t)$}\\
1- \eps^{-1} (|x|-R - C_0 (\tau-t)) & \mbox{if $0\leq |x|- (R+C_0 (\tau-t)) \leq \eps$}\\
0 &\mbox{otherwise,}
\end{array}\right.
$$
where $C_0$ is the constant appearing in \eqref{e:dafermos_5.3.5}.
Note that:
\begin{itemize}
\item $0\leq \psi\leq 1$;
\item $\psi (t,x) = 0$ if $t\geq \tau$ or $|x|\geq \eps+R + C (\tau-t)$;
\item $\partial_t \psi = -\eps^{-1}$ on $B_R (0)\times ]\tau-\eps, \tau[$;
\item $|\nabla_x \psi|\leq - C_0^{-1} \partial_t \psi$.
\end{itemize}
Combining these pieces of information with \eqref{e:dafermos_5.3.5},
from \eqref{e:dafermos_5.3.10} we easily conclude
\begin{equation}
\frac{1}{\eps} \int_{\tau-\eps}^\tau \int_{|x|\leq R}
h\, dx\, dt \;\leq\; \int_0^\tau \int_{|x|\leq R+\eps + C_0 (\tau-t)} |\nabla U||Z|\,
dx\, dt
\end{equation}
Recalling \eqref{e:dafermos_5.3.5_bis} and the Lipschitz regularity
of $U$ we conclude
\begin{equation}
\frac{1}{\eps} \int_{\tau-\eps}^\tau \int_{|x|\leq R}
h\, dx\, dt \;\leq\; C \int_0^\tau \int_{|x|\leq R+\eps + C_0 (\tau-t)} h\, dx\, dt\, .
\end{equation}
Finally, letting $\eps\downarrow 0$ and using the fact that $h$ is integrable, we conclude
\begin{equation}\label{e:dafermos_5.3.14}
\int_{|x|\leq R} h (x,\tau)\, dx \;\leq\;
C \int_0^\tau \int_{|x|\leq R+ C_0 (\tau-t)} h(x,t)\, dx\, dt
\qquad \mbox{for a.e. $\tau$.}
\end{equation}
Note, moreover, that the set of measure zero where \eqref{e:dafermos_5.3.14}
fails can be chosen independently of $R$. Therefore, having fixed any $s<T$,
we infer
\begin{equation}\label{e:one_more}
\int_{|x|\leq R + C_0 (s-\tau)} h (x,\tau)\, dx \;\leq\;
C \int_0^\tau \int_{|x|\leq R+ C_0 (s-t)} h(x,t)\, dx\, dt\quad
\mbox{for a.e. $\tau\in [0,s]$.}
\end{equation}
If we set
$$
g(\tau):=\int_{|x|\leq R + C_0 (s-\tau)} h (x, \tau)\, dx\, ,
$$
then \eqref{e:one_more} becomes the Gronwall's inequality 
$g(\tau)\leq C \int_0^\tau g(t)\, dt$,
which leads to the conclusion $g\equiv 0$. By the arbitrariness of $R>0$ and $s<T$ we
conclude that $h\equiv 0$ on $[0,T]\times \R^n$. Recalling \eqref{e:other_side},
we infer $\nu_{x,t} = \delta_{U(x,t)}$ for a.e. $(x,t)$, which is the desired conclusion.
\end{proof}


\end{document}